\documentclass{article}

 \twocolumn 

\usepackage{times}
\usepackage{amsmath}
\usepackage{amssymb}
\usepackage{latexsym}

\newcommand{\pro}{P}
\newcommand{\pr}[1]{P\!\left( #1 \right)}

\newcommand{\lpr}[1]{\underline{P}\!\left( #1 \right)}

\newcommand{\upr}[1]{\overline{P}\!\left( #1 \right)}

\newcommand{\exoo}{E}
\newcommand{\exo}[1]{E\!\left[ #1 \right]}
\newcommand{\ex}[2]{E_{#1}\!\left[ #2 \right]}

\newcommand{\lex}[1]{\underline{E}\!\left[ #1 \right]}
\newcommand{\uexo}{\overline{E}}
\newcommand{\uex}[1]{\overline{E}\!\left[ #1 \right]}

\newcommand{\credalo}{K}

\newtheorem{Theorem}{Theorem}

\newtheorem{Example}{Example}

\title{Concentration Inequalities and 
       Laws of Large Numbers under Epistemic Irrelevance}

\author{Fabio Gagliardi Cozman \\
        Escola Polit\'ecnica, 
        Universidade de S\~ao Paulo -
        S\~ao Paulo, SP - Brazil \\
        {\tt fgcozman@usp.br}}

\begin{document}

\maketitle

\noindent
{\bf NOTE: this is a revised version of a paper
presented at ISIPTA 2009. This is not a major revision;
the purpose was to correct some problems with the original paper.
The main point is that Example 1 has been changed (the original example 
was flawed). Also, the definition of $B_i$ and some inequalities
in the proof of limits in Theorem 4 have been corrected.}

\vspace*{3ex}

\begin{abstract}
This paper presents concentration inequalities and laws 
of large numbers under weak assumptions of irrelevance,
expressed through lower and upper expectations. The results 
are variants and extensions of De Cooman and Miranda's recent 
inequalities and laws of large numbers. The proofs indicate 
connections between concepts
of irrelevance for lower/upper expectations and the standard 
theory of martingales.
\end{abstract}

\section{Introduction}
\label{section:Introduction}

This paper investigates concentration inequalities 
and laws of large numbers under weak assumptions of 
``irrelevance'' that are expressed using lower and 
upper expectations. The starting point is the assumption
that, given bounded variables $X_1, \dots, X_n$, we have:
%
%
\begin{equation}
\label{equation:AssumptionEpistemicIrrelevance}
\parbox{6.5cm}{
\centerline{for each $i \in [2,n]$, variables $X_1,\dots,X_{i-1}$}

\centerline{are epistemically irrelevant to $X_i$.}}
\end{equation}
Epistemic irrelevance of variables $X_1,\dots,X_{i-1}$ to 
$X_i$ obtains when \cite[Def. 9.2.1]{Walley91}
\begin{equation}
\label{equation:EpistemicIrrelevance}
\uex{f(X_i)|A(X_{1:i-1})} = \uex{f(X_i)}
\end{equation}
for any bounded function $f$ of $X_i$ and any nonempty event 
$A(X_{1:i-1})$ defined by variables $X_{1:i-1}$, where
the functional $\uexo$ is 
an {\em upper expectation} (Section \ref{section:Basics}). 
Here and in the remainder of the paper we simplify
notation by using $X_{1:i}$ for $X_1,\dots,X_i$.

A judgement of epistemic irrelevance can be interpreted 
as a relaxed judgement of stochastic independence, perhaps
motivated by a robustness analysis or by disagreements amongst
a set of decision makers. Alternatively, one might 
consider epistemic irrelevance as {\em the} appropriate
concept of independence when expectations are not known precisely.

De Cooman and Miranda have recently proven a number of inequalities 
and laws of large numbers that also deal with judgements 
of irrelevance expressed through lower/upper expectations
\cite{Cooman2008}. 
De Cooman and Miranda's weak law of large numbers implies 
that, given Assumption (\ref{equation:AssumptionEpistemicIrrelevance}),
 for any $\epsilon > 0$, 
\[
\lpr{ \underline{\mu}_n - \epsilon \leq
      \frac{\sum_{i=1}^n X_i}{n} \leq
       \overline{\mu}_n + \epsilon } \geq  
             1 - 2e^{-\frac{n\epsilon^2/4}{(\max_i B_i)^2}},
\]
where $B_i$ is such that $\sup X_i - \inf X_i \leq B_i$, and
\[
 \underline{\mu}_n \doteq \frac{\sum_{i=1}^n \lex{X_i}}{n}, 
\qquad
 \overline{\mu}_n  \doteq \frac{\sum_{i=1}^n \uex{X_i}}{n}.
\]
Moreover, De Cooman and Miranda's results and Assumption 
(\ref{equation:AssumptionEpistemicIrrelevance}) imply 
a two-part strong law of large numbers: for any $\epsilon > 0$, there 
is $N \in \mathbf{N}_{+}$ such that for any $N' \in \mathbf{N}_{+}$,
\[
\upr{ \exists n \in [N, N+N'] :  
   \frac{\sum_{i=1}^n X_i}{n} \geq \overline{\mu} + \epsilon} < \epsilon,
\]
\[
\upr{ \exists n \in [N, N+N'] :  
   \frac{\sum_{i=1}^n X_i}{n} \leq \underline{\mu} - \epsilon} < \epsilon.
\]
This law of large numbers corresponds to a finitary version 
of the usual strong law of large numbers \cite{Dubins74}; the
focus on a finitary law is justified by the fact that De Cooman
and Miranda do not assume countable additivity. 
If countable additivity 
holds, the finitary strong law of large numbers implies
convergence of empirical means with probability one
\cite[Sec. 5.3]{Cooman2008}.

To obtain their results, De Cooman and Miranda assume, following
Walley's theory of lower previsions, that all variables are bounded,
and that conglomerability (and consequently disintegrability) holds. 
These assumptions are discussed in more detail later.

The present paper derives laws of large numbers by exploiting
concentration and martingale inequalities that are adapted to
the setting of lower/upper expectations. These results use either
Assumption (\ref{equation:AssumptionEpistemicIrrelevance}) or
the weaker assumption that, for each $i \in [2,n]$ and any nonempty 
event $A(X_{1:i-1})$,
\begin{equation}
\label{equation:AssumptionWeakIrrelevance}
\parbox{6.5cm}{
\centerline{$\lex{X_i|A(X_{1:i-1})} = \lex{X_i}$} 
 
\centerline{and} 

\centerline{$\uex{X_i|A(X_{1:i-1})} = \uex{X_i}$.}}
\end{equation}
Several results for {\em bounded} variables presented in this paper are 
basically implied by De Cooman and Miranda's work. Regarding bounded
variables our contribution lies in offering tighter inequalities 
and alternative proof techniques that are more closely related
to established methods in standard probability theory (in 
particular, close to Hoeffding's and Azuma's inequalities).
In Section \ref{section:Unbounded} we offer more significant
contributions as we lift the assumption of boundedness for variables,
and use martingale theory to prove laws of large numbers under 
elementwise disintegrability. 

\section{Expectations, disintegrability, and zero probabilities}
\label{section:Basics}

In this section we present notation and terminology.
Throughout the paper we assume that an expectation functional
$\exoo$  maps bounded variables into real numbers, and satisfies: \\
(1) if $\alpha \leq X \leq \beta$, then $\alpha \leq \exo{X} \leq \beta$; \\
(2) $\exo{X+Y} = \exo{X} + \exo{Y}$; \\
where $X, Y$ are bounded variables and $\alpha, \beta$ are real numbers
(inequalities are understood pointwise).

From such an expectation functional, a {\em finitely additive}
probability measure $\pro$ is induced by $\pr{A} \doteq \exo{A}$ for
any event $A$; note that $A$ denotes both the event and its indicator
function.\footnote{A probability measure defined on a field completely
characterizes an expectation functional on bounded functions that
are measurable with respect to the field and vice-versa 
\cite[Theorem 3.2.2]{Walley91}.}
%

Given a set of expectation functionals, 
the lower and upper expectations of variable $X$ are respectively
\[
\lex{X} = \inf \exo{X}, \qquad \uex{X} = \sup \exo{X}.
\]
Lower and upper probabilities are defined similarly
using indicator functions. Given an event $A$, a conditional
expectation functional is constrained by 
$\exo{X|A} \pr{A} = \exo{XA}$. If we have a set of expectation 
functionals, then a set of conditional expectation functionals 
given an event $A$ is produced by elementwise conditioning on 
event $A$ (that is, each expectation functional is conditioned on $A$).

\subsection{Disintegrability and factorization}
\label{subsection:Disintegrability}

We will employ an assumption of {\em disintegrability} 
in our proofs; namely,
\begin{equation}
\label{equation:Conglomerability}
\uex{W} \leq \uex{\uex{W|Z}} 
\end{equation}
for any $W \geq 0$, $Z \geq 0$ of interest, where $W$ and $Z$ 
may stand for sets of (non-negative) variables. Note that 
disintegrability can fail for a single finitely additive
probability measure over an infinite space \cite{Finetti74,Dubins75};
that is, there is a finitely additive probability measure $\pro$ 
such that
\[
\ex{\pro}{W} > \ex{\pro}{\ex{\pro}{W|Z}}.
\]
One way to obtain disintegrability is to restrict attention to 
simple variables; that is, variables that take on finitely many 
distinct values. In particular, indicator functions are simple variables;
hence simple variables suffice to express convergence of 
relative frequencies, and our results apply then.

Another way to obtain disintegrability for every probability
measure $\pro$ is to adopt countable additivity \cite{Ash99}.
That is, assume that if
\[ 
A_1 \supset A_2 \supset \dots
\]
is a countable sequence of events, then
\begin{equation}
\label{equation:CountableAdditivity}
\cap_i A_i = \emptyset \quad \mbox{ implies } \quad
\lim_{n \rightarrow \infty} \upr{A_n} = 0.
\end{equation}
This assumption says that if $\cap_i A_i = \emptyset$, then
$\lim_{n \rightarrow \infty} \pr{A_n} = 0$ for every possible 
probability measure.

A third way to obtain disintegrability is simply to impose it.
One may consider disintegrability a ``rationality'' requirement. 
\begin{itemize}
\item The theories of coherent behavior by Heath and Sudderty~\cite{Heath78} 
and by Lane and Sudderth~\cite{Lane85} follow this path by axiomatizing 
the {\em strategic} measures of Dubins and Savage~\cite{Dubins65}, 
and thus prescribing probability measures 
that disintegrate appropriately along some predefined partitions.
This would be sufficient for our purposes, but there are limitations
in the approach as summarized by Kadane et al \cite{Kadane96ASA}.
The disintegrability of strategic measures has actually been used to
prove various laws of large numbers in a finitely additive setting
\cite{Karandikar82}. 
\item Another scheme that imposes disintegrability is Walley's 
theory of lower previsions; in that theory,
Expression (\ref{equation:Conglomerability}) is a consequence
of axioms for ``coherent'' behavior. This is the path
adopted by De Cooman and Miranda, who consequently have 
Expression (\ref{equation:Conglomerability}) at their disposal.
\end{itemize}

When disintegrability holds, recursive application of
Expression (\ref{equation:Conglomerability}) yields:
if $f_i(X_i) \geq 0$ for $i \in \{1,\dots,n\}$, then
\begin{eqnarray*}
\uex{\prod_{i=1}^n f_i(X_i)} & & \\
& & \hspace{-2cm} \leq 
   \uex{ \dots \uex{\uex{ \prod_{i=1}^n f_i(X_i) | X_{1:n-1}}|X_{1:n-2}}\dots};
\end{eqnarray*}
Assumption (\ref{equation:AssumptionEpistemicIrrelevance}) 
then implies an inequality we use later: for bounded and nonnegative functions,
\begin{equation}
\label{equation:Factorization}
\uex{\prod_{i=1}^n f_i(X_i)} \leq \prod_{i=1}^n \uex{f_i(X_i)}.
\end{equation}

\subsection{Zero probabilities, full conditional measures
and weak irrelevance}
\label{subsection:Zero}

It should be noted that the definition of epistemic irrelevance 
(Expression (\ref{equation:EpistemicIrrelevance})) does not contain 
any clause concerning zero probabilities. 
Indeed, Walley's theory of lower previsions follows 
de Finetti in adopting {\em full conditional measures}, and in 
this setting Expression (\ref{equation:EpistemicIrrelevance}) 
can be imposed without concerns about zero probabilities. 
Recall that a full conditional measure 
$\pro: \mathcal{B} \times (\mathcal{B}\backslash\emptyset) 
	\rightarrow \Re$, where $\mathcal{B}$ is a Boolean 
algebra, is a set-function that for every nonempty event $C$
satisfies~\cite{Dubins75,Krauss68}: \\
(1) $\pr{C|C} = 1$;   \\
(2) $\pr{A|C} \geq 0$ for all $A$;  \\
(3) $\pr{A \cup B|C}=\pr{A|C}+\pr{B|C}$ for all 
	disjoint $A$ and $B$; \\
(4) $\pr{A \cap B|C} = \pr{A|B \cap C} \pr{B|C}$ for all $A$ 
	and $B$ such that $B \cap C \neq \emptyset$.  

Full conditional measures are not adopted in the usual
Kolmogorovian theory, and if countable additivity is 
adopted and co´nditioning is defined through Radon-Nykodym derivatives, 
it may be impossible to satisfy the axioms for full conditional measures 
\cite{Seidenfeld2001Synthese,Seidenfeld2001AP}. Thus there are are
some differences between epistemic irrelevance (at least as
defined by Walley) and the usual Kolmogorovian set-up, besides 
the obvious set-valued/point-valued distinction.

Suppose that one wishes to deal with sets of probability measures
and associated lower/upper expectations,
but chooses to adopt the Kolmogorovian set-up for each measure.
That is, each measure satisfies countable additivity and thus
disintegrability, and conditioning is left undefined when the
conditioning event has probability zero.
It might seem reasonable to amend
Expression (\ref{equation:EpistemicIrrelevance}) as follows:
\begin{eqnarray}
\label{equation:TooWeak}
\uex{f(X_i)|A(X_{1:i-1})} & = & \uex{f(X_i)} \\
& &  \mbox{if $\lpr{A(X_{1:i-1})}>0$}. \nonumber
\end{eqnarray}
This condition is a natural for theories that do not define
conditioning on events of lower probability zero, such as 
Giron and Rios' theory \cite{Giron80}.
Alas, this weaker condition is really too weak to produce
laws of large numbers, as the following example shows.

\begin{Example}
Consider binary variables $X_1, X_2, \dots$ (values $0$ and $1$).
Define events
$A_0 \doteq \{X_1=0,X_2=0,\dots\}$ and
$A_1 \doteq \{X_1=1,X_2=1,\dots\}$.
Consider a convex and closed set $\credalo$ of joint distributions
for these variables, built as the convex hull of three distributions,
$\pro_1$, $\pro_2$ and $\pro_3$, as follows.

Distribution $\pro_1$ simply assigns probability one to $A_1$.
Distribution $\pro_2$ assigns probability $\delta$ to $A_0$
and probability $1-\delta$ to $A_1$, for some $\delta \in (0,1)$.
Distribution $\pro_3$ is the product of
identical marginals: for any integer $n>0$,
$\pro_3(X_1=x_1,\dots,X_n=x_n) = \prod_{i=1}^n \pro_3(X_i=x_i)$,
where $\pro_3(X_i=1) = 1-\delta$.

For the convex hull of $\pro_1$, $\pro_2$ and $\pro_3$, 
Expression (\ref{equation:TooWeak}) is satisfied.
This conclusion is reached by analyzing each distribution in turn.
For distribution $\pro_1$, we have $\pro_1(X_1=1)=1$ and 
for any $i>1$ we have $\pro_1(X_i=1|A(X_{1:i-1}))=1$ whenever
$\lpr{A(X_{1:i-1})}>0$. 
Note that for any event $A(X_{1:i-1})$:
if $A_1 \in A$, then $\pro_1(A) = 1$;
if $A_1 \not\subseteq A$, then $\pro_1(A) = 0$.
For distribution $\pro_2$, 
$\pro_2(X_i=1)=1-\delta$ for any $i>0$. 
Additionally, for any event $A(X_{1:i-1})$ 
we have $\pro_2(X_i=1|A)$ either equal 
to $1-\delta$ or $1$ whenever $\lpr{A}>0$.
[If $A_1 \not\subseteq A$, then $\lpr{A} = 0$ (due to $\pro_1$).
So suppose $A_1 \subseteq A$:
If $A_0 \subseteq A$, then $\pro_2(X_i=1|A) = 1-\delta$;
if $A_0 \not\subseteq A$, then $\pro_2(X_i=1|A) = 1$.]
For distribution $\pro_3$, we have
$\pro_3(X_i=1)=1-\delta$ and for any $i>1$ we have
$\pro_3(X_i=1|A)=1-\delta$
for any nonempty event $A(X_{1:i-1})$.
In short, for all probability measures in the credal set
we have $\pr{X_i=1} \in [1-\delta,1]$ and
$\pr{X_i=1|A(X_{1:i-1})} \in [1-\delta,1]$ whenever $\lpr{A(X_{1:i-1})}>0$..

The weak law of larger numbers fails because,
for any $\epsilon \in (0,1-\delta)$,
\[
\lim_{n \rightarrow \infty}
\lpr{ 
\underline{\mu}_n - \epsilon \leq
\frac{\sum_{i=1}^n X_i}{n} \leq
\overline{\mu}_n + \epsilon } = 1-\delta.
\]
This follows from the fact that, for any integer $n>0$, 
we have $\pro_1\left(\sum_{i=1}^n X_i/n = 1\right) = 1$ and 
$\pro_2\left( \sum_{i=1}^n X_i/n = 1 \right) = 1-\pro_2(A_0) = 1-\delta$,
and for any $\epsilon>0$ (due to standard weak law of large numbers),
\[
\lim_{n \rightarrow \infty}
\pro_3\left((1-\delta)-\epsilon < \sum_{i=1}^n X_i/n  
                                < (1-\delta)+\epsilon\right) = 1.
\]
\end{Example}

We might thus consider 
an alternative to Expression~(\ref{equation:TooWeak}):
\begin{eqnarray}
\label{equation:Regular}
\uex{f(X_i)|A(X_{1:i-1})} & = & \uex{f(X_i)} \\
& & \mbox{if $\upr{A(X_{1:i-1})}>0$}. \nonumber
\end{eqnarray}
The concept of irrelevance conveyed by 
Expression (\ref{equation:Regular}) does lead to
Expression (\ref{equation:Factorization}). To see this,
note that for nonnegative $X$ and $Y$, we have
\begin{eqnarray*}
\uex{XY} & \leq & \sup_{\pro} \ex{\pro}{\uex{XY|Y}} \\
         & = & \sup_{\pro} \ex{\pro}{ A \uex{XY|Y} + A^c \uex{XY|Y} },
\end{eqnarray*}
using disintegrability and defining $A$ 
as the set of all values of $Y$ such that $\upr{A^c} = 0$.
Hence $\pr{A^c} = 0$ for every $\pro$ and using
Expression (\ref{equation:Regular}):
\begin{eqnarray*}
\uex{XY} &  \leq & \sup_{\pro} \ex{\pro}{A Y \uex{X|Y}} \\
         & =     & \sup_{\pro} \ex{\pro}{A Y \uex{X}} \\
         & =     & \sup_{\pro} \ex{\pro}{A Y} \uex{X} \\
         & =     & \uex{X} \sup_{\pro} \ex{\pro}{Y}  \\
         & =     & \uex{X} \uex{Y}.
\end{eqnarray*}

[As a digression, note that one might {\em define} 
conditional expectations as 
$\lex{X|A} = \inf_{\pro : \pr{A}>0} \ex{\pro}{X|A}$ and 
$\uex{X|A} = \sup_{\pro : \pr{A}>0} \ex{\pro}{X|A}$.
This form of conditioning has been advocated by 
several authors \cite{Weichselberger2000,Weichselberger2001},
and it is quite similar to Walley's concept of regular 
extension \cite[Ap. J]{Walley91}. For such a form of
conditioning, Expression (\ref{equation:Regular}) seems to
be the natural definition of irrelevance.]

In short, more than one combination of definitions and assumptions
lead to the results presented in the remainder of this paper. 
For instance, Expression (\ref{equation:Factorization}) obtains when
Assumption (\ref{equation:AssumptionEpistemicIrrelevance}) holds 
{\em and} disintegrability holds 
(because all variables are simple, {\em or} because 
countable additivity is assumed, {\em or} because disintegrability 
is imposed).
Alternatively, Expression (\ref{equation:Factorization}) obtains when
Expression (\ref{equation:Regular}) holds for any $i \in [2,n]$, any
bounded function $f$ of $X_i$, and any event $A(X_{1:i-1})$, and 
additionally disintegrability holds.

Similar remarks concerning zero probabilities can be directed 
at Assumption (\ref{equation:AssumptionWeakIrrelevance}). 
We say that {\em weak irrelevance} obtains when either one of:
\begin{itemize}
\item For any $i \in [2,n]$ and any nonempty event $A(X_{1:i-1})$, 
\begin{center}
$\lex{X_i|A(X_{1:i-1})} = \lex{X_i}$ \\ and \\
$\uex{X_i|A(X_{1:i-1})} = \uex{X_i}$
\end{center}
[this is Assumption (\ref{equation:AssumptionWeakIrrelevance}),
and it requires full conditional measures].
\item For any $i \in [2,n]$ and any event $A(X_{1:i-1})$, 
\begin{center}
$\lex{X_i|A(X_{1:i-1})} = \lex{X_i}$ if $\upr{A(X_{1:i-1})}>0$
\\ and \\
$\uex{X_i|A(X_{1:i-1})} = \uex{X_i}$ if $\upr{A(X_{1:i-1})}>0$.
\end{center}
\end{itemize}

\section{Bounded variables}
\label{section:Bounded}

Take variables $X_1,\dots,X_n$ such that $\sup X_i - \inf X_i \leq B_i$
and define 
\[
\gamma_n \doteq \sum_{i=1}^n B_i^2 > 0.
\] 
We start by deriving two concentration inequalities.

\subsection{Concentration inequalities}

The following inequality is a counterpart of Hoeffding
inequality \cite{Devroye96,Hoeffding63} in the context of
lower/upper expectations; it is slightly tighter than similar
inequalities by De Cooman and Miranda \cite{Cooman2008}.
It is interesting to note that  
the proof is remarkably similar to the proof of the original 
Hoeffding inequality. 

\begin{Theorem}
\label{theorem:HoeffdingExtended}
If bounded variables $X_1, \dots, X_n$ satisfy 
Expression (\ref{equation:Factorization}), then
if $\gamma_n>0$,
\[
\upr{\sum_{i=1}^n (X_i - \uex{X_i}) \geq \epsilon} 
     \leq e^{-2\epsilon^2/\gamma_n},
\]
\[
\upr{\sum_{i=1}^n (X_i - \lex{X_i}) \leq -\epsilon} 
     \leq e^{-2\epsilon^2/\gamma_n}.
\]
\end{Theorem}
{\em Proof.}
By Markov inequality, if $X \geq 0$, then for any $\epsilon>0$
we have $\pr{X \geq \epsilon} \leq \exo{X}/\epsilon$. 
Consequently, for $s>0$, any variable $X$ satisfies
\[
\upr{X \geq \epsilon} = \upr{e^{sX} \geq e^{s\epsilon}} \leq
e^{-s\epsilon} \uex{\exp(sX)}.
\]
Using this inequality and Expression (\ref{equation:Factorization}):
\begin{eqnarray*}
\upr{\sum_{i=1}^n (X_i-\uex{X_i}) \geq \epsilon } & & \\
& & \hspace*{-3.5cm} \leq 
e^{-s\epsilon}
 \uex{\exp\left(\sum_{i=1}^n s(X_i-\uex{X_i})\right)} \\
& & \hspace*{-3.5cm} \leq 
e^{-s\epsilon}
     \prod_{i=1}^n \uex{\exp\left(s(X_i-\uex{X_i})\right)}.
\end{eqnarray*}
We now use Hoeffding's result (Expression (\ref{equation:HoeffdingResult}))
that if variable $X$ satisfies $a \leq X \leq b$ 
and $\exo{X} \leq 0$, then
$\exo{\exp(sX)} \leq \exp(s^2(b-a)^2/8)$ for any $s>0$.
Thus for any $\pro$, 
$\ex{\pro}{\exp(s(X_i-\uex{X_i}))} \leq \exp(s^2 B_i^2/8)$,
and then 
$\uex{\exp\left(s(X_i-\uex{X_i})\right)} \leq \exp(s^2B_i^2/8)$.
Consequently,
\[
\upr{\sum_{i=1}^n (X_i-\uex{X_i}) \geq \epsilon } \leq
 e^{-s\epsilon}
 e^{s^2 \gamma_n /8} \leq  e^{-2\epsilon^2/\gamma_n},
\]
where the last inequality is obtained by taking
$s = 4 \epsilon/\gamma_n$. This proves the first inequality
in the theorem; the second inequality is proved by 
taking $\upr{\sum_{i=1}^n ((-X_i)-\uex{-X_i}) \geq \epsilon }$
and noting that $\lex{X_i} = - \uex{-X_i}$.
$\Box$

We now move to weak irrelevance and obtain an analogue
of Azuma's inequality \cite{Chung2006,Devroye91}. 
It is again interesting to note that the proof is 
remarkably similar to the proof of the
 original Azuma inequality.
De Cooman and Miranda \cite[Sec. 4.1]{Cooman2008} show that their
inequalities are valid under weak irrelevance; the next inequality
is slightly tighter than theirs.

\begin{Theorem}
\label{theorem:AzumaExtended}
If bounded variables $X_1, \dots, X_n$ satisfy weak irrelevance
and disintegrability (Expression (\ref{equation:Conglomerability}))
holds, then if $\gamma_n>0$,
\[
\upr{ \sum_{i=1}^n (X_i - \uex{X_i}) \geq \epsilon }
     \leq e^{-2\epsilon^2/\gamma_n},
\]
\[
\upr{ \sum_{i=1}^n (X_i - \lex{X_i}) \leq -\epsilon }
     \leq e^{-2\epsilon^2/\gamma_n}.
\]
\end{Theorem}
{\em Proof.} 
Using both Markov's inequality (as in the proof of 
Theorem \ref{theorem:HoeffdingExtended}) and 
disintegrability, for any $s>0$ we get
\begin{eqnarray*}
\upr{\sum_{i=1}^n (X_i-\uex{X_i}) \geq \epsilon } & & \\
& & \hspace*{-4.5cm} \leq 
e^{-s\epsilon}
 \uex{\exp\left(\sum_{i=1}^n s(X_i-\uex{X_i})\right)} \\
& & \hspace*{-4.5cm} \leq 
e^{-s\epsilon}
 \uex{\uex{\exp\!\!\left(\sum_{i=1}^n s(X_i-\uex{X_i})\right) 
                                         \mid X_{1:n-1}}} \\
& & \hspace*{-4.5cm} \leq 
 e^{-s\epsilon}
 \uex{\exp\!\!\left(\sum_{i=1}^{n-1} s(X_i-\uex{X_i})\right) h(X_{1:n-1})},
\end{eqnarray*}
where
\[
h(X_{1:n-1}) = 
 \uex{ \exp\!\!\left(s(X_n-\uex{X_n})\right) \mid X_{1:n-1} }.
\]
Due to weak irrelevance,
\[
\ex{\pro}{X_n|X_{1:n-1}} \leq \uex{X_n|X_{1:n-1}} = \uex{X_n};
\]
consequently, for any $\pro$,
\[
\ex{\pro}{X_n-\uex{X_n}|X_{1:n-1}} \leq 0.
\]
We now use Hoeffding's result (Expression (\ref{equation:HoeffdingResult}))
that if variable $X$ satisfies $a \leq X \leq b$ 
and $\exo{X} \leq 0$, then
$\exo{\exp(sX)} \leq \exp(s^2(b-a)^2/8)$ for any $s>0$.
Thus for any $\pro$ we have
\[
\ex{\pro}{\exp\left(s(X_n-\uex{X_n})\right)|X_{1:n-1}} \leq \exp(s^2 B_n^2/8)
\]
and then $h(X_{1:n-1}) \leq \exp(s^2B_n^2/8)$. Thus
\begin{eqnarray*}
\upr{\sum_{i=1}^n (X_i-\uex{X_i}) \geq \epsilon } & & \\
& & \hspace*{-4cm} \leq 
e^{-s\epsilon}
 \uex{\exp\left(\sum_{i=1}^n s(X_i-\uex{X_i})\right)} \\
& & \hspace*{-4cm} \leq 
 e^{-s\epsilon}
  \uex{\exp\!\!\left(\sum_{i=1}^{n-1} s(X_i-\uex{X_i})\right) 
  \exp(s^2B_n^2/8)} \\
& & \hspace*{-4cm} \leq 
 e^{-s\epsilon}   \exp(s^2B_n^2/8)
  \uex{\exp\!\!\left(\sum_{i=1}^{n-1} s(X_i-\uex{X_i})\right)}. 
\end{eqnarray*}
These inequalities can be iterated to produce:
\[
\upr{\sum_{i=1}^n (X_i-\uex{X_i}) \geq \epsilon } \leq 
 e^{-s\epsilon}   \exp\left(s^2 \sum_{i=1}^n B_i^2/8\right).
\]
Finally, by taking $s = 4\epsilon/\gamma_n$,
\[
\upr{\sum_{i=1}^n (X_i-\uex{X_i}) \geq \epsilon } \leq 
         e^{-2\epsilon^2/\gamma_n}.
\]
The second inequality in the theorem is proved by
noting that weak irrelevance of $X_1,\dots,X_n$ implies
weak irrelevance of $-X_1,\dots,-X_n$ 
(as $\lex{X_i} = - \uex{-X_i}$), and then by
taking $\upr{\sum_{i=1}^n ((-X_i)-\uex{-X_i}) \geq \epsilon }$.
$\Box$

\subsection{Laws of large numbers}

Theorem \ref{theorem:HoeffdingExtended} leads to simple
proofs of laws of large numbers already stated by 
De Cooman and Miranda \cite{Cooman2008}. To start, take
Assumption (\ref{equation:AssumptionEpistemicIrrelevance}).
Using subadditivity of upper probability 
and Theorem~\ref{theorem:HoeffdingExtended},
\[
\upr{\!\!\!\left(\sum_{i=1}^n X_i \geq n\overline{\mu}_n + \epsilon\!\right) 
     \!\cup\!
     \left(\sum_{i=1}^n X_i \leq n\underline{\mu}-\epsilon\!\right)\!\!\!}
     \!\leq\!2e^{-\frac{2\epsilon^2}{\gamma_n}},
\]
where as before, 
$\underline{\mu}_n \doteq (1/n)\sum_{i=1}^n \lex{X_i}$ and 
$\overline{\mu}_n  \doteq (1/n)\sum_{i=1}^n \uex{X_i}$.
By noting that $\lpr{A} = 1 - \upr{A^c}$ for any event $A$, 
by including the endpoints of relevant inequalities,
and by using $n \epsilon$ instead of $\epsilon$:
\begin{eqnarray*}
\lpr{ \underline{\mu} - \epsilon 
              \leq \frac{\sum_{i=1}^n X_i}{n} \leq
      \overline{\mu} + \epsilon } & \geq & \\
\lpr{ \underline{\mu} - \epsilon 
              < \frac{\sum_{i=1}^n X_i}{n} <
      \overline{\mu} + \epsilon }  & \geq & 
             1 - 2e^{-\frac{2n\epsilon^2}{B^2}},
\end{eqnarray*}
where we define $B \doteq \max_i B_i$. By taking limits,
we obtain a weak law of large numbers:
\[
\lim_{n \rightarrow \infty}
\lpr{ \underline{\mu}_n - \epsilon 
             < \frac{\sum_{i=1}^n X_i}{n} <
      \overline{\mu}_n + \epsilon } = 1.
\]
An analogue of De Cooman and Miranda's finitary strong law of 
large numbers can be deduced as well from the previous inequalities,
as follows. Here and in the remainder of the paper,
$n$, $N$ and $N'$ denote positive integers. For all
$\epsilon > 0$, $N > 0$ and $N'>0$, 
\begin{eqnarray*}
\upr{ \exists n \in [N, N+N'] :  
   \frac{\sum_{i=1}^n X_i}{n} \geq \overline{\mu} + \epsilon}
 & & \\
& & \hspace*{-4.5cm}
\leq 
\sum_{n=N}^{N+N'} \upr{
   \frac{\sum_{i=1}^n X_i}{n} \geq \overline{\mu} + \epsilon} \\
& & \hspace*{-4.5cm}
\leq 
\sum_{n=N}^{N+N'} e^{- 2 n \epsilon^2/B^2} \\
& & \hspace*{-4.5cm}
=
\left( e^{-2N\epsilon^2/B^2} \right) 
  \sum_{n=0}^{N'} e^{-2n\epsilon^2/B^2} \\
& & \hspace*{-4.5cm}
=
\left( e^{-2N\epsilon^2/B^2} \right) 
   \frac{1-e^{2(N'+1)\epsilon^2/B^2}}
        {1 - e^{-2\epsilon^2/B^2}} \\
& & \hspace*{-4.5cm}
 <  \frac{e^{-2N\epsilon^2/B^2}}{1 - e^{-2\epsilon^2/B^2}}.
\end{eqnarray*}
Consequently,
\[
\upr{ \exists n \in [N, N+N'] :  
   \frac{\sum_{i=1}^n X_i}{n} \geq \overline{\mu} + \epsilon} < \epsilon,
\]
provided that $N$ is a positive integer such that
\[
N > - (B^2/(2\epsilon^2)) \ln \epsilon(1-e^{-2\epsilon^2/B^2}).
\]
An analogous argument leads to
\[
\upr{ \exists n \in [N, N+N'] :  
   \frac{\sum_{i=1}^n X_i}{n} \leq \underline{\mu} - \epsilon} < \epsilon.
\]
By superadditivity of upper probability, we obtain a perhaps
more intuitive statement of the strong law of large numbers:
for all $\epsilon>0$, there is $N$ such that for any $N'$,
\[
\lpr{\!\forall n\!\in\![N, N\!+\!N']\!:\!%
\underline{\mu}_n\!-\!\epsilon\!<\!\frac{\sum_{i=1}^n X_i}{n}\!<\!\overline{\mu}_n\!+\!\epsilon\!}\!>\!1-2\epsilon,
\]
thus reproducing De Cooman and Miranda's strong laws.

We now present a pair of weak/strong laws of large numbers 
under weak irrelevance. De Cooman and Miranda prove a similar 
pair of laws by resorting to their previous results on 
{\em forward irrelevant natural extensions} \cite[Sec. 4.1]{Cooman2008}. 
The proof offered now is perhaps more direct, using our
analogue of Azuma's inequality.

\begin{Theorem}
\label{theorem:LLNbounded}
If bounded variables $X_1, \dots, X_n$ satisfy weak irrelevance
and Expression (\ref{equation:Conglomerability}) holds, then
for any $\epsilon > 0$, 
\[
\lpr{ \underline{\mu}_n - \epsilon 
              < \frac{\sum_{i=1}^n X_i}{n} < 
       \overline{\mu}_n + \epsilon } \geq  
             1 - 2e^{-2n\epsilon^2/B^2},
\]
and there is $N$ such that for any $N'$,
\[
\lpr{\!\forall n\!\in\![N,N\!+\!N']\!:\!%
\underline{\mu}_n\!-\!\epsilon\!<\!\frac{\sum_{i=1}^n\!X_i}{n}\!<\!\overline{\mu}_n\!+\!\epsilon\!}\!\!>\!\!1-2\epsilon.
\]
\end{Theorem}
{\em Proof.}
Using subadditivity of upper probability 
and Theorem \ref{theorem:AzumaExtended}, and defining
again $B \doteq \max_i B_i$,
\[
\upr{\!\!\!\left(\sum_{i=1}^n X_i \geq n\overline{\mu}_n + \epsilon\!\right) 
     \!\!\cup\!\!
     \left(\sum_{i=1}^n X_i \leq n\underline{\mu}-\epsilon\!\right)\!\!\!}\!%
\leq\!2e^{-\frac{2n\epsilon^2}{B^2}},
\]
and we obtain the first expression in the theorem. 
To produce the second inequality (strong law), note:
\begin{eqnarray*}
\upr{ \exists n \in [N, N+N'] :  
 \frac{\sum_{i=1}^n X_i}{n} \geq \overline{\mu} + \epsilon}
 & & \\
& & \hspace*{-4.5cm}
\leq 
\sum_{n=N}^{N+N'} \upr{
 \frac{\sum_{i=1}^n X_i}{n} \geq \overline{\mu} + \epsilon} \\
& & \hspace*{-4.5cm}
\leq 
\sum_{n=N}^{N+N'} e^{-2 n \epsilon^2/B^2} \\
& & \hspace*{-4.5cm}
 <  \frac{e^{-2N\epsilon^2/B^2}}{1 - e^{-2\epsilon^2/B^2}}.
\end{eqnarray*}
Again,
\[
\upr{ \exists n \in [N, N+N'] :  
 \frac{\sum_{i=1}^n X_i}{n} \geq \overline{\mu} + \epsilon} < \epsilon
\]
provided that $N$ is a positive integer such that
\[
N > - (B^2/(2\epsilon^2)) \ln \epsilon(1-e^{-2\epsilon^2/B^2}).
\]
This is ``half'' of the second expression in the theorem;
the other ``half'' is proved analogously. 
$\Box$

The theorem easily implies the following concise weak
law of large numbers, by taking limits:
\[
\lim_{n \rightarrow \infty}
\lpr{ \underline{\mu}_n - \epsilon 
              < \frac{\sum_{i=1}^n X_i}{n} <
      \overline{\mu}_n + \epsilon } = 1.
\]

\section{Laws of large numbers without boundedness}
\label{section:Unbounded}

We now consider variables without bounds in their ranges
under the assumption of weak irrelevance; the resulting
laws of large numbers are the main contribution 
of the paper. We will assume in this section that countable 
additivity holds (Expression (\ref{equation:CountableAdditivity})).
This assumption of countable addivity implies disintegrability; 
that is, $\ex{\pro}{W} = \ex{\pro}{\ex{\pro}{W|Z}}$ 
for any $\pro$, $W$ and $Z$.
Thus our setup is close to the standard (Kolmogorovian) one, 
where any expectation functional is a linear monotone and monotonically 
convergent functional that can be expressed
through Lebesgue integration. We only depart from the Kolmogorovian
tradition in explicitly letting a {\em set} of such functionals to be
permissible given a set of assessments.


We will use a sequence of variables $\{Y_n\}$ defined as follows:
\[
Y_n \doteq \sum_{i=1}^n X_i - \ex{\pro}{X_i|X_{1:i-1}}.
\]
The key observation is that $Y_n$ is a function of all variables 
$X_{1:n}$ such that
\begin{eqnarray*}
\ex{\pro}{Y_n | X_{1:n-1} } & = & 
 \left( \sum_{i=1}^{n-1} X_i - \ex{\pro}{X_i|X_{1:i-1}} \right) + \\
& & \ex{\pro}{ X_n\!-\!\ex{\pro}{X_n|X_{1:n-1}} | X_{1:n-1} } \\
& = & Y_{n-1} + \\
& & \ex{\pro}{ X_n | X_{1:n-1} } - \ex{\pro}{X_n|X_{1:n-1}} \\
& = & Y_{n-1};
\end{eqnarray*}
so, $\{Y_n\}$ is a {\em martingale} with respect to $\pro$.
Thus,
\begin{eqnarray*}
\ex{\pro}{(Y_n-Y_{n-1})^2|X_{1:n-1}} & & \\
& & \hspace*{-4.5cm} = 
\ex{\pro}{Y_n^2|X_{1:n-1}} - 2 \ex{\pro}{Y_{n-1}Y_n|X_{1:n-1}} + Y_{n-1}^2 \\
& & \hspace*{-4.5cm} = 
\ex{\pro}{Y_n^2|X_{1:n-1}} - 2 Y_{n-1} \ex{\pro}{Y_n|X_{1:n-1}} + Y_{n-1}^2 \\
& & \hspace*{-4.5cm} = 
\ex{\pro}{Y_n^2|X_{1:n-1}} - 2 Y_{n-1} Y_{n-1} + Y_{n-1}^2 \\
& & \hspace*{-4.5cm} = 
\ex{\pro}{Y_n^2|X_{1:n-1}} - Y_{n-1}^2.
\end{eqnarray*}
And by taking expectations on both sides and noting
that $Y_i-Y_{i-1} = X_i - \ex{\pro}{X_i|X_{1:i-1}}$, we get
\[
\ex{\pro}{Y^2_n} = \ex{\pro}{(X_n-\ex{\pro}{X_n|X_{1:n-1}})^2} 
                  + \ex{\pro}{Y^2_{n-1}}.
\]
Iterating this expression, we obtain:
\begin{equation}
\label{equation:MartingaleIteration}
\ex{\pro}{Y^2_n} = \sum_{i=1}^n \ex{\pro}{(X_i-\ex{\pro}{X_i|X_{1:i-1}})^2}.
\end{equation}

With these preliminaries, we have: 
\begin{Theorem}
\label{theorem:LLNunbounded}
Assume countable additivity.
If variables $X_1, \dots, X_n$ satisfy weak irrelevance, and
$\lex{X_i}$ and $\uex{X_i}$ are finite quantities such 
that $\uex{X_i}-\lex{X_i} \leq \delta$, and
the variance of any $X_i$ is no larger than a finite quantity $\sigma^2$, 
then for any $\epsilon > 0$, 
\[
\lpr{ \underline{\mu}_n - \epsilon 
              < \frac{\sum_{i=1}^n X_i}{n} <
       \overline{\mu}_n + \epsilon } \geq 1 - \frac{\sigma^2+\delta^2}
                                                   {\epsilon^2 n},
\]
and there is $N>0$ such that for any $N'>0$,
\[
\lpr{\!\forall n\!\in\![N,N\!+\!N']\!:\!%
\underline{\mu}_n\!-\!\epsilon\!<\!\frac{\sum_{i=1}^n\!X_i}{n}\!<\!\overline{\mu}_n\!+\!\epsilon\!}\!\!>\!\!1-2\epsilon.
\]
Consequently,
\[
\forall \epsilon>0: \;
\lim_{n \rightarrow \infty}
\lpr{ \underline{\mu}_n - \epsilon 
             < \frac{\sum_{i=1}^n X_i}{n} <
      \overline{\mu}_n + \epsilon } = 1,
\]
\[
\lpr{ \lim\sup_{n \rightarrow \infty} 
 \left(  \frac{\sum_{i=1}^n X_i}{n} - \overline{\mu}_n \right) \leq 0 } = 1,
\]
\[
\lpr{ \lim\inf_{n \rightarrow \infty} 
 \left(  \frac{\sum_{i=1}^n X_i}{n} - \underline{\mu}_n \right) \geq 0 } = 1.
\]
\end{Theorem}
{\em Proof.}
For a fixed $\pro$ and for all $\epsilon>0$,
\begin{eqnarray*}
\pr{ \underline{\mu}_n - \epsilon 
              < \frac{\sum_{i=1}^n X_i}{n} <
       \overline{\mu}_n + \epsilon } & & \\
& & \hspace*{-5.5cm} = 
\pr{\sum_{i=1}^n \lex{X_i} - \epsilon n <
        \sum_{i=1}^n X_i <
    \sum_{i=1}^n \uex{X_i} + \epsilon n } \\
& & \hspace*{-5.5cm} \geq
\pro\left(\sum_{i=1}^n \ex{\pro}{X_i|X_{1:i-1}} - \epsilon n <
        \sum_{i=1}^n X_i \right. \\
& & \hspace*{-2.5cm} 
   \left. < \sum_{i=1}^n \ex{\pro}{X_i|X_{1:i-1}} + \epsilon n \right) \\
& & \hspace*{-5cm} \mbox{ (using weak irrelevance)} \\
& & \hspace*{-5.5cm} =
\pr{-\epsilon < \frac{\sum_{i=1}^n X_i - \ex{\pro}{X_i|X_{1:i-1}}}{n}
              < \epsilon} \\
& & \hspace*{-5.5cm} =
\pr{ -\epsilon < Y_n/n < \epsilon } \\
& & \hspace*{-5.5cm} = 
\pr{|Y_n/n| < \epsilon}.
\end{eqnarray*}
Applying Chebyshev's inequality and
Expression (\ref{equation:MartingaleIteration}),
\begin{eqnarray*}
\pr{|Y_n/n| \geq \epsilon} & \leq & 
       \frac{\ex{\pro}{Y_n^2}}{\epsilon^2 n^2} \\
& = & 
 \frac{\sum_{i=1}^n\!\ex{\pro}{(X_i\!-\!\ex{\pro}{X_i|X_{1:i-1}})^2}}
                {\epsilon^2 n^2}.
\end{eqnarray*}
Now write $(X_i-\ex{\pro}{X_i|X_{1:i-1}})^2$ as
\[
\left((X_i-\ex{\pro}{X_i})+(\ex{\pro}{X_i}-\ex{\pro}{X_i|X_{1:i-1}})\right)^2,
\] 
and then:
\begin{eqnarray*}
\sum_{i=1}^n \ex{\pro}{(X_i-\ex{\pro}{X_i|X_{1:i-1}})^2} & & \\
& & \hspace*{-5.5cm} =
\sum_{i=1}^n \ex{\pro}{(X_i-\ex{\pro}{X_i})^2}  \\
& & \hspace*{-5cm} 
             + 2 \ex{\pro}{(X_i-\ex{\pro}{X_i})
                         (\ex{\pro}{X_i}- \ex{\pro}{X_i|X_{1:i-1}})} \\
& & \hspace*{-4.5cm} 
             + \ex{\pro}{(\ex{\pro}{X_i}- \ex{\pro}{X_i|X_{1:i-1}})^2} \\
& & \hspace*{-5.5cm} \leq
\sum_{i=1}^n \sigma^2  + \delta^2 \\
& & \hspace*{-5cm} 
             + 2 (\ex{\pro}{X_i}- \ex{\pro}{X_i|X_{1:i-1}})
                    \ex{\pro}{X_i-\ex{\pro}{X_i}} \\
& & \hspace*{-5.5cm} =
\sum_{i=1}^n \sigma^2  + \delta^2.
\end{eqnarray*}
Hence
\begin{equation}
\label{equation:SigmaDelta}
\sum_{i=1}^n \ex{\pro}{(X_i-\ex{\pro}{X_i|X_{1:i-1}})^2} \leq 
n (\sigma^2 + \delta^2),
\end{equation}
and combining these inequalities, we obtain:
\[
\pr{|Y_n/n| \geq \epsilon} \leq \frac{\sigma^2+\delta^2}{\epsilon^2 n}, 
\]
and then
\[
\pr{ \underline{\mu}_n - \epsilon 
              < \frac{\sum_{i=1}^n X_i}{n} <
       \overline{\mu}_n + \epsilon } \geq 1 - \frac{\sigma^2+\delta^2}
                                                   {\epsilon^2 n}
\]
for any $\pro$, as desired. By taking the limit as $n$ grows without
bound, we obtain
\[
\lim_{n \rightarrow \infty}
\lpr{ \underline{\mu}_n - \epsilon 
              < \frac{\sum_{i=1}^n X_i}{n} <
      \overline{\mu}_n + \epsilon } = 1.
\]


The proof of the strong law of large numbers uses the same
strategy, but replaces the appeal to Chebyshev's inequality
by an appeal to the Kolmogorov-Hajek-Renyi inequality
(described in the Appendix), following the proof of the strong
law of large numbers by Whittle \cite[Thm. 14.2.3]{Whittle92}.
So, for a fixed $\pro$ and for all $\epsilon>0$, we proceed
as previously to obtain:
\begin{eqnarray*}
\pr{\!\forall n\!\in\![N,N\!+\!N']: 
\underline{\mu}_n\!-\!\epsilon\!< \frac{\sum_{i=1}^n\!X_i}{n}
   <\!\overline{\mu}_n\!+\!\epsilon\!} & & \\
& & \hspace*{-5.5cm} \geq
\pr{\!\forall n\!\in\![N,N\!+\!N']: 
-\epsilon < \frac{Y_n}{n} < \epsilon} \\
& & \hspace*{-5.5cm} =
\pr{\!\forall n\!\in\![N,N\!+\!N']: |Y_n/n| < \epsilon }.
\end{eqnarray*}
As $\{Y_N, Y_{N+1}, \dots, Y_{N+N'}\}$ forms a martingale,
we use the Kolmogorov-Hajek-Renyi inequality to produce:
\begin{eqnarray*}
\pr{\!\forall n\!\in\![N,N\!+\!N']: |Y_n/n| < \epsilon } & & \\
& & \hspace*{-4cm} \geq
1 - \frac{ \sum_{i=1}^N \ex{\pro}{(X_i - \ex{\pro}{X_i|X_{1:i-1}})^2} }
         {\epsilon^2 N^2} \\
& & \hspace*{-3cm} 
- \sum_{i=N+1}^{N+N'} \frac{ \ex{\pro}{(X_i - \ex{\pro}{X_i|X_{1:i-1}})^2}}
                           {\epsilon^2 i^2} \\
& & \hspace*{-4cm} \geq 
1 - \frac{ \sigma^2 + \delta^2 }{\epsilon^2 N} 
  - \sum_{i=N+1}^{N+N'} \frac{\sigma^2 + \delta^2}{\epsilon^2i^2} \\
& & \hspace*{-3cm} \mbox{ (using Expression (\ref{equation:SigmaDelta}))} \\
& & \hspace*{-4cm} \geq 
1 - \frac{ \sigma^2 + \delta^2 }{\epsilon^2 N} 
  - \sum_{i=N+1}^{\infty} \frac{\sigma^2 + \delta^2}{\epsilon^2i^2} \\
& & \hspace*{-4cm} \geq 
1 - \frac{ \sigma^2 + \delta^2 }{\epsilon^2}
  \left(\frac{1}{N}+\int_{N}^\infty 1/i^2 di\right) \\
& & \hspace*{-4cm} =
1 - \frac{ \sigma^2 + \delta^2 }{\epsilon^2}
  \left(\frac{1}{N}+\frac{1}{N}\right) \\
& & \hspace*{-4cm} =
1 - 2\frac{ \sigma^2 + \delta^2 }{\epsilon^2 N}. 
\end{eqnarray*}
Consequently, for integer
$N > (\sigma^2+\delta^2)/\epsilon^3$, we obtain the desired
inequality
\[
\lpr{\!\forall n\!\in\![N,N\!+\!N']\!:\!%
\underline{\mu}_n\!-\!\epsilon\!< \frac{\sum_{i=1}^n\!X_i}{n}
  <\!\overline{\mu}_n\!+\!\epsilon\!}\!\!>\!\!1-2\epsilon.
\]
The proof of the Kolmogorov-Hajek-Renyi can be
extended to an infinite intersection of (decreasing) events
expressed as $\{\forall j \geq 1 : |X_j| < \epsilon_j\}$; thus
\begin{eqnarray*}
\forall \epsilon>0 : \forall \delta > 0 : \exists N > 0: & & \\
& & \hspace*{-4cm}
\pr{\forall m \geq N : \frac{\sum_{i=1}^m X_i - \uex{X_i}}{m} <
                      \epsilon} \geq 1-\delta,
\end{eqnarray*}
and this is equivalent to:
\[
\forall \epsilon>0:
\lim_{N \rightarrow \infty} \pr{\!\forall m\!\geq\!N :
\frac{\sum_{i=1}^m X_i - \uex{X_i}}{m} < \epsilon\!} = 1.
\]
As the events in these probability values form an increasing sequence,
we have, for all $\epsilon > 0$,
\[
\pr{\exists N > 0: \forall m \geq N :
     \frac{\sum_{i=1}^m X_i - \uex{X_i}}{m} < \epsilon} = 1. 
\]
Now this is equivalent to $\forall k > 0 : \pr{A_k} = 1$,
where $A_k = \{ \exists N > 0: \forall m \geq N: 
     (1/m) \sum_{i=1}^m X_i - \uex{X_i} > 1/k \}$,
and because $\pr{\cup_{k > 0} \neg A_k} \leq \sum_{k>0}\pr{\neg A_k}=0$,
we have $\pr{\forall k>0: A_k} = 1$, so
\[
\pr{\!\!\forall k\!>\!0\!:\exists N\!>\!0\!:\forall m\!\geq\!N\!:\!%
\frac{\sum_{i=1}^m X_i - \uex{X_i}}{m} < \epsilon\!}=1.
\]
This is exactly the desired expression
\[
\lpr{ \lim\sup_{n \rightarrow \infty} 
      \left(  \frac{\sum_{i=1}^n X_i}{n} 
                    - \overline{\mu}_n \right) \leq 0 } = 1.
\]
A similar argument proves the last inequality in the theorem, starting
from: 
\begin{eqnarray*}
\forall \epsilon>0 : \forall \delta > 0 : \exists N > 0: & & \\
& & \hspace*{-4cm}
\pr{\forall m \geq N : \frac{\sum_{i=1}^m X_i - \lex{X_i}}{m} 
			> - \epsilon} \geq 1-\delta.
\end{eqnarray*}
$\Box$

\section{Discussion}

The concentration inequalities and laws of large numbers
proved in this paper assume rather weak conditions of epistemic
irrelevance. When compared to usual laws of large numbers, both premises
and consequences are weaker: expectations are not assumed 
precisely known, and convergence is interval-valued.

Theorems \ref{theorem:HoeffdingExtended} and 
\ref{theorem:AzumaExtended} and their ensuing 
laws of large numbers are implied by De Cooman and Miranda's 
seminal work \cite{Cooman2008} (and their results generalize
several previous efforts \cite{Epstein2003}). Actually, De Cooman
and Miranda start from a weaker condition of
{\em forward factorization} that is implied both by Assumption
(\ref{equation:AssumptionEpistemicIrrelevance}) and 
weak irrelevance.
The possible advantage of our proof techniques for these
two theorems is that they
are rather close to well-known methods in standard probability
theory, such as Hoeffding's inequality (it should be noted that 
De Cooman and Miranda already indicate the similarity between
their inequalities and Hoeffding's).

The most significant results of the paper employ weak irrelevance
to produce concentration inequalities (Theorem \ref{theorem:AzumaExtended})
and laws of large numbers (Theorems \ref{theorem:LLNbounded} and
\ref{theorem:LLNunbounded}). The latter theorem is possibly the
most valuable contribution. The strategy for most proofs is to translate
assumptions of weak irrelevance into facts regarding martingales, and
to adapt  results for martingales to this setting. This strategy keeps 
the proof relatively short and close to well-known results in probability 
theory. The connection between lower/upper expectations and the theory 
of martingales seems rather natural \cite{Cooman2008AI,Shafer2001}, 
but the relationship between epistemic irrelevance and martingales
does not appear to have been explored in depth so far.
We note that the basic constraint defining martingales (that is, 
$\exo{Y_n|X_{1:n-1}}=Y_{n-1}$) is preserved by convex
combination of mixtures; therefore, the study of martingales 
seems appropriate when one deals with convex sets of probability
measures --- certainly it seems less contorted than the analysis
through stochastic independence, as stochastic independence is
{\em not} preserved by convex combination.

The proofs presented in this paper need assumptions of disintegrability
that can be easily satisfied if countable additivity is adopted.
It is an open question whether similar results can be proven without
disintegrability, particularly when one deals with unbounded variables.



\section*{Acknowledgements}

Thanks to Gert de Cooman and Enrique Miranda for 
their generous suggestions regarding
content and presentation. Thanks to a reviewer who indicated
the work by Sadrolhefazi and Fine on ergodic 
properties \cite{Sadrolhefazi94}, a topic to be pursued
in the future.

The author is partially supported by CNPq.

\appendix

\section{Two auxiliary inequalities}

The following inequality is a simple extension of 
a basic result by Hoeffding \cite{Devroye96,Hoeffding63}: 
If variable $X$ satisfies $a \leq X \leq b$ 
and $\exo{X} \leq 0$, then for any $s>0$, 
\begin{equation}
\label{equation:HoeffdingResult} 
\exo{\exp(sX)} \leq \exp(s^2 (b-a)^2/8).
\end{equation}
First, the inequality is clearly valid if $a = b$,
or if $a=0$, or if $b<0$.
From now on, suppose $b \geq 0 > a$.
By convexity of the exponential function,
\[
\exp(sx) \leq \frac{x-a}{b-a}e^{sb} + \frac{b-x}{b-a}e^{sa}
\quad \mbox{ for $x \in [a,b]$}.
\]
Given monotonicity of expectations and $\exo{X} \leq 0$,
\[
\exo{\exp(sX)} \leq \frac{b}{b-a}e^{sa} - \frac{a}{b-a}e^{sb}
\doteq \exp(\phi(s(b-a)))
\]
for $\phi(u) = -pu+\log(1-p+pe^u)$ with $p=-a/(b-a)$ (and note
that $p \in (0,1]$ in the situation under consideration).
Given that $\phi(0)=\phi'(0)=0$ and $\phi''(u) \leq 1/4$ for $u>0$
(as the maximum of $\phi''(u)$ is $1/4$, attained at $e^u = (1-p)/p$),
we can use Taylor's theorem as follows. For some
$v \in (0,u)$, 
$\phi(u) = \phi(0)+u\phi'(0)+(u^2/2)\phi''(v) \leq (u^2/8)$ 
and consequently $\phi(s(b-a)) \leq s^2(b-a)^2/8$.
By putting together these inequalities, we obtain Expression
(\ref{equation:HoeffdingResult}).

We now review the Kolmogorov-Hajek-Renyi inequality, almost exactly
as proved by Whittle \cite{Whittle92}; this is presented just
to indicate the role of (elementwise) disintegrability in the derivation.
Let $\{X_i\}$ be a martingale with $X_0 = 0$, and let $\{\epsilon_i\}$
be a sequence $0=\epsilon_0 \leq \epsilon_1 \leq \dots$; the  
inequality is
\[
\pr{\forall j \in [1,n] : |X_j| < \epsilon_j} \geq 
    1 - \sum_{i=1}^n \frac{\exo{(X_i-X_{i-1})^2}}{\epsilon_i^2}.
\]
To prove this inequality, define
$A_n \doteq \{ \forall j \in [1,n] : |X_j| < \epsilon_j\}$.
Using $\xi_i = X_i - X_{i-1}$, and again denoting an event
and its indicator function by the same symbol, we have
\begin{eqnarray*}
\pr{A_n} & = & \ex{\pro}{A_n} \;\; = \;\;
              \ex{\pro}{A_{n-1} \{|X_n| < \epsilon_n\}} \\
         & \geq & \ex{\pro}{A_{n-1} (1 - X_n^2/\epsilon_n^2)}  \\
         &   &   \mbox{ (as $\{|X| < \epsilon\} \geq 1-X^2/\epsilon^2$) } \\
         & = & \ex{\pro}{A_{n-1} (1 - (X_{n-1}^2 + \xi_n^2)/\epsilon_n^2)} \\
         &   &   \mbox{ (by the martingale property) } \\
         & \geq & \ex{\pro}{A_{n-2} (1 - X_{n-1}^2/\epsilon_{n-1}^2)} -
                     \ex{\pro}{\xi_n^2/\epsilon_n^2} \\
         &   &   \mbox{ (as $\epsilon_{n-1} \leq \epsilon_n$ and} \\
         &   &   \mbox{ 
$\{ |X| < \epsilon \} (1 - X^2/\epsilon^2) \geq (1-X^2/\epsilon^2)$)}.
\end{eqnarray*}
Iteration of the last inequality yields the result.
Note that it was necessary to apply disintegrability of $\pro$ 
when applying the martingale property
(that is, elementwise disintegrability is used).

\end{document}